\documentclass[12pt]{amsart}
\usepackage{amssymb,amsmath,amsfonts,latexsym,setspace}
\usepackage{bm}
\newcommand{\newsection}[1]{\setcounter{equation}{0} \section{#1}}
\newcommand{\bea}{\begin{eqnarray}}
\newcommand{\eea}{\end{eqnarray}}

\newcommand{\clb}{\mathcal{B}}

\newcommand{\cle}{\mathcal{E}}
\newcommand{\clf}{\mathcal{F}}

\newcommand{\clh}{\mathcal{H}}
\newcommand{\clk}{\mathcal{K}}
\newcommand{\cll}{\mathcal{L}}
\newcommand{\clm}{\mathcal{M}}

\newcommand{\clq}{\mathcal{Q}}
\newcommand{\clr}{\mathcal{R}}
\newcommand{\cls}{\mathcal{S}}

\newcommand{\raro}{\rightarrow}

\def \qed {\hfill \vrule height6pt width 6pt depth 0pt}
\def\textmatrix#1&#2\\#3&#4\\{\bigl({#1 \atop #3}\ {#2 \atop #4}\bigr)}
\def\dispmatrix#1&#2\\#3&#4\\{\left({#1 \atop #3}\ {#2 \atop #4}\right)}
\newcommand{\be}{\begin{equation}}
\newcommand{\ee}{\end{equation}}
\newcommand{\ben}{\begin{eqnarray*}}
\newcommand{\een}{\end{eqnarray*}}

\newcommand{\NI}{\noindent}

\newcommand{\bi}{\begin{itemize}}
\newcommand{\ei}{\end{itemize}}

\newtheorem{Theorem}{\sc Theorem}[section]
\newtheorem{Lemma}[Theorem]{\sc Lemma}
\newtheorem{Proposition}[Theorem]{\sc Proposition}
\newtheorem{Corollary}[Theorem]{\sc Corollary}
\newtheorem{Definition}[Theorem]{\sc Definition}
\newtheorem{Example}[Theorem]{\sc Example}
\newtheorem{Remark}[Theorem]{\sc Remark}
\newtheorem{Note}[Theorem]{\sc Note}
\newtheorem{Question}{\sc Question}
\newtheorem{ass}[Theorem]{\sc Assumption}
\newcommand{\bt}{\begin{Theorem}}
\def\beginlem{\begin{Lemma}}
\def\beginprop{\begin{Proposition}}
\def\begincor{\begin{Corollary}}
\def\begindef{\begin{Definition}}
\def\beginexamp{\begin{Example}}
\def\beginrem{\begin{Remark}}
\def\beginq{\begin{Question}}
\def\beginass{\begin{ass}}
\def\beginnote{\begin{Note}}
\newcommand{\et}{\end{Theorem}}
\def\endlem{\end{Lemma}}
\def\endprop{\end{Proposition}}
\def\endcor{\end{Corollary}}
\def\enddef{\end{Definition}}
\def\endexamp{\end{Example}}
\def\endrem{\end{Remark}}
\def\endq{\end{Question}}
\def\endass{\end{ass}}
\def\endnote{\end{Note}}

\begin{document}

\title[Submodules of the Hardy module over polydisc]{Submodules of the Hardy module over polydisc}

\author[Jaydeb Sarkar]{Jaydeb Sarkar}
\address{Indian Statistical Institute, Statistics and Mathematics Unit, 8th Mile, Mysore Road, RVCE Post, Bangalore, 560059, India}
\email{jay@isibang.ac.in, jaydeb@gmail.com}

\dedicatory{Dedicated to Ronald G. Douglas on the occasion of his
75th birthday}



\begin{abstract}
We say that a submodule $\cls$ of $H^2(\mathbb{D}^n)$ ($n >1$) is
co-doubly commuting if the quotient module $H^2(\mathbb{D}^n)/\cls$
is doubly commuting. We show that a co-doubly commuting submodule of
$H^2(\mathbb{D}^n)$ is essentially doubly commuting if and only if
the corresponding one variable inner functions are finite Blaschke
products or that $n = 2$. In particular, a co-doubly commuting
submodule $\cls$ of $H^2(\mathbb{D}^n)$ is essentially doubly
commuting if and only if $n = 2$ or that $\cls$ is of finite
co-dimension. We obtain an explicit representation of the
Beurling-Lax-Halmos inner functions for those submodules of
$H^2_{H^2(\mathbb{D}^{n-1})}(\mathbb{D})$ which are co-doubly
commuting submodules of $H^2(\mathbb{D}^n)$. Finally, we prove that
a pair of co-doubly commuting submodules of $H^2(\mathbb{D}^n)$ are
unitarily equivalent if and only if they are equal.
\end{abstract}

\maketitle

\newsection{Introduction}
Let $\{T_1, \ldots, T_n\}$ be a set of $n$ commuting bounded linear
operators on a separable Hilbert space $\clh$. Then we can turn the
$n$-tuple $(T_1, \ldots, T_n)$ on $\clh$ into a Hilbert module
\cite{DP} $\clh$ over $\mathbb{C}[\bm{z}] := \mathbb{C}[z_1, \ldots,
z_n]$, the ring of polynomials, as follows:
\[\mathbb{C}[\bm{z}] \times \clh \raro \clh, \quad \quad (p, h)
\mapsto p(T_1, \ldots, T_n)h,\]for all $p \in \mathbb{C}[\bm{z}]$
and $h \in \clh$. The \textit{module multiplication} operators by
the coordinate functions on $\clh$ are defined by $M_{z_i} h =
z_i(T_1, \ldots, T_n) h = T_i h$, for all $h \in \clh$ and $i = 1,
\ldots, n$. Therefore, a Hilbert module is uniquely determined by
the underlying commuting operators via the module multiplication
operators by the coordinate functions and vice versa. Let $\cls,
\clq \subseteq \clh$ be closed subspaces of $\clh$. Then $\cls$
($\clq$) is said to be a \textit{submodule} (\textit{quotient
module}) of $\clh$ if $M_{z_i} \cls \subseteq \cls$ ($M^*_{z_i} \clq
\subseteq \clq$) for all $i = 1, \ldots, n$. Note that a closed
subspace $\clq$ is a quotient module of $\clh$ if and only if
$\clq^{\perp} \cong \clh/\clq $ is submodule of $\clh$.

The \textit{Hardy module} $H^2(\mathbb{D}^n)$ over the polydisc is
the Hardy space $H^2(\mathbb{D}^n)$ (cf. \cite{FS} and \cite{R}),
the closure of $\mathbb{C}[\bm{z}]$ in $L^2(\mathbb{T}^n)$, with the
standard multiplication operators by the coordinate functions $z_i$
($1 \leq i \leq n$) on $H^2(\mathbb{D}^n)$ as the module maps.

The module multiplication operators on a submodule $\cls$ and a
quotient module $\clq$ of a Hilbert module $\clh$ are given by the
restrictions $(R_{z_1}, \ldots, R_{z_n})$ and the compressions
$(C_{z_1}, \ldots, C_{z_n})$ of the module multiplications of
$\clh$, respectively. That is,
\[R_{z_i} = M_{z_i}|_{\cls} \quad \mbox{and} \quad C_{z_i} = P_{\clq}
M_{z_i}|_{\clq},\]for all $i = 1, \ldots, n$. Here, for a given
closed subspace $\clm$ of a Hilbert space $\clk$, we denote the
orthogonal projection of $\clk$ onto $\clm$ by $P_{\clm}$.

A quotient module $\clq$ of a Hilbert module $\clh$ over
$\mathbb{C}[\bm{z}]$ ($n\geq 2$) is said to be \textit{doubly
commuting quotient module} if
\[C_{z_i} C_{z_j}^* = C_{z_j}^* C_{z_i},\] for all $1 \leq i < j
\leq n$. Also a submodule $\cls$ of $\clh$ is said to be
\textit{co-doubly commuting submodule} of $\clh$ if $\clh/\cls$ is
doubly commuting quotient module of $\clh$ (see \cite{JS}).

Finally, we recall that a Hilbert module $\clh$ over
$\mathbb{C}[\bm{z}]$ is said to be \textit{essentially doubly
commuting} if the cross-commutators
\[ [M_{z_i}^*, M_{z_j}] \in \clk(\clh),\]for all $1 \leq i < j \leq n$,
where $\clk(\clh)$ is the ideal of all compact operators on $\clh$.
We say that $\clh$ is \textit{essentially normal} if $[M_{z_i}^*,
M_{z_j}] \in \clk(\clh)$ for all $1 \leq i, j \leq n$.

Natural examples of essentially normal Hilbert modules are the
Drury-Arveson module $H^2_n$, the Hardy module $H^2(\mathbb{B}^n)$
and the Bergman module $L^2_a(\mathbb{B}^n)$ over the unit ball
$\mathbb{B}^n$ (cf. \cite{D}, \cite{A}, \cite{A-trans}). On the
other hand, the Hardy module $H^2(\mathbb{D}^n)$ over $\mathbb{D}^n$
with $n \geq 2$ is not essentially normal. However, a simple
calculation reveals that $H^2(\mathbb{D}^n)$ is doubly commuting
and, in particular, essentially doubly commuting. Therefore, a
natural approach to measure a submodule (quotient module) of the
Hardy module $H^2(\mathbb{D}^n)$ from being "small" is to consider
the cross commutators $[R_{z_i}^*, R_{z_j}]$ ($[C_{z_i}^*,
C_{z_j}]$) for all $1 \leq i < j \leq n$ instead of all possible
commutators.

Before proceeding further, let us recall the Beurling-Lax-Halmos
theorem concerning submodules of vector-valued Hardy modules over
$\mathbb{D}$ (cf. \cite{NF}).

\NI Given a separable Hilbert space $\cle$ we shall denote by
$H^2_{\cle}(\mathbb{D})$ the $\cle$-valued Hardy module (see
\cite{NF}). Note that by virtue of the unitary module map $U :
H^2_{\cle}(\mathbb{D}) \raro H^2(\mathbb{D}) \otimes \cle$ defined
by \[z^m \eta \mapsto z^m \otimes \eta, \quad (\eta \in \cle, m \in
\mathbb{N})\] we can identify the vector-valued Hardy module
$H^2_{\cle}(\mathbb{D})$ with $H^2(\mathbb{D}) \otimes \cle$.

\begin{Theorem}\label{BLHthm}\textit{(Beurling-Lax-Halmos)}
Let $\cle$ be a Hilbert space and $\cls$ be a non-trivial closed
subspace of the Hardy module $H^2_{\cle}(\mathbb{D})$. Then $\cls$
is a submodule of $H^2_{\cle}(\mathbb{D})$ if and only if
\[\cls = \Theta H^2_{\clf}(\mathbb{D}),\]where $\Theta \in
H^\infty_{\cll(\clf, \cle)}(\mathbb{D})$ is an inner function and
$\clf$ is an adequate Hilbert space with dimension less than or
equal to the dimension of $\cle$. Moreover, $\Theta$ is unique up to
a unitary constant right factor, that is, if $\cls = \tilde{\Theta}
H^2_{\tilde{\clf}}(\mathbb{D})$ for some Hilbert space
$\tilde{\clf}$ and inner function $\tilde{\Theta} \in
H^\infty_{\cll(\tilde{\clf}, \cle)}(\mathbb{D})$, then $\Theta =
\tilde{\Theta} W$ where $W$ is a unitary operator in $\cll(\clf,
\tilde{\clf})$.
\end{Theorem}

Now we formulate some general problems concerning submodules of the
Hardy module $H^2(\mathbb{D}^n)$ ($n \geq 2$).

\begin{Question}\label{Q1}\textsf{(Essentially doubly commuting submodules)}
How to characterize essentially doubly commuting
submodules of the Hardy modules $H^2(\mathbb{D}^n)$?
\end{Question}

Let $\cls \neq \{0\}$ be a closed subspace of
$H^2_{H^2(\mathbb{D}^{n-1})}(\mathbb{D})$. By Beurling-Lax-Halmos
theorem, Theorem \ref{BLHthm}, that $\cls$ is a submodule of
$H^2_{H^2(\mathbb{D}^{n-1})}(\mathbb{D})$ if and only if $\cls =
\Theta H^2_{\cle_*}(\mathbb{D})$, for some closed subspace $\cle_*
\subseteq H^2(\mathbb{D}^{n-1})$ and inner function $\Theta \in
H^\infty_{\cll(\cle_*, H^2(\mathbb{D}^{n-1}))}(\mathbb{D})$.

\begin{Question}\label{Q2}\textsf{(Beurling-Lax-Halmos representations)} For
which closed subspace $\cle_* \subseteq H^2(\mathbb{D}^{n-1})$ and
inner function $\Theta \in H^\infty_{\cll(\cle_*,
H^2(\mathbb{D}^{n-1}))}(\mathbb{D})$ the submodule $\Theta
H^2_{\cle_*}(\mathbb{D})$ of
$H^2_{H^2(\mathbb{D}^{n-1})}(\mathbb{D})$, realized as a subspace of
$H^2(\mathbb{D}^n)$, is a submodule of $H^2(\mathbb{D}^n)$?
\end{Question}

Let $\clh$ be a Hilbert module over $\mathbb{C}[\bm{z}]$. Denote by
$\clr(\clh)$ the set of all non-unitarily equivalent submodules of
$\clh$, that is, if $\cls_1, \cls_2 \in \clr(\clh)$ and that $\cls_1
\cong \cls_2$ then $\cls_1 = \cls_2$.

\begin{Question}\label{Q3}\textsf{(Rigidity of submodules)}
Determine $\clr(H^2(\mathbb{D}^n))$.
\end{Question}

The aim of the present paper is to analyze and answer the above
questions for the class of co-doubly commuting submodules of
$H^2(\mathbb{D}^n)$. We obtain an explicit description of the cross
commutators of co-doubly commuting submodules of
$H^2(\mathbb{D}^n)$. As an applications, we prove that the cross
commutators of a co-doubly commuting submodule $\cls$ of
$H^2(\mathbb{D}^n)$ are compact, that is, $\cls$ is essentially
doubly commuting, if and only if $n = 2$ or that $\cls$ is of finite
co-dimension.  We would like to point out that a submodule of finite
co-dimension is necessarily essentially doubly commuting. Therefore,
the issue of essential doubly commutativity of co-doubly commuting
submodules of $H^2(\mathbb{D}^n)$ yields a rigidity type result: if
$\cls$ is of infinite co-dimension co-doubly commuting submodules of
$H^2(\mathbb{D}^n)$ and $\cls$ is essentially doubly commuting, then
$n = 2$ (the base case). Our earlier classification results are also
used to prove a Beurling-Lax-Halmos type theorem for the class of
co-doubly commuting submodules of $H^2(\mathbb{D}^n)$. We also
discuss the rigidity phenomenon of such submodules.

Note also that most of the results of the present paper, concerning
doubly commuting quotient modules and co-doubly commuting
submodules, restricted to the base case $n= 2$ are known. However,
the proofs are new even in the case $n = 2$. Moreover, as we have
pointed out above, the difference between the base case $n=2$ and
the higher variables case $n>2$ is more curious in the study of
essentially doubly commuting submodules of $H^2(\mathbb{D}^n)$ (see
Corollaries \ref{EN-2}, \ref{<3} and \ref{n=2}).

We now summarize the contents of this paper. In Section 2 we
investigate the essential doubly commutativity problem for the class
of co-doubly commuting submodules of $H^2(\mathbb{D}^n)$ and
conclude that for $n \geq 3$, except for the finite co-dimension
case, none of the co-doubly commuting submodules of
$H^2(\mathbb{D}^n)$ are essentially doubly commuting. In Sections 3
and 4 we answer Questions \ref{Q2} and \ref{Q3} for the class of
co-doubly commuting submodules of $H^2(\mathbb{D}^n)$, respectively.
We conclude in Section 5 with some remarks and discussion on the
problem of essentially doubly commutativity of Hilbert modules.

\newsection{Cross commutators of submodules}

In a recent paper \cite{JS} we completely classify the class of
doubly commuting quotient modules and co-doubly commuting submodules
of the Hardy module $H^2(\mathbb{D}^n)$, where $n \geq 2$ (see
\cite{INS} for the case $n = 2$).

\begin{Theorem}\label{S-proj}
Let $\clq$ be a quotient module of $H^2(\mathbb{D}^n)$ and $n \geq
2$. Then $\clq$ is doubly commuting quotient module of
$H^2(\mathbb{D}^n)$ if and only if
\[\clq = \clq_{\Theta_1} \otimes \cdots \otimes \clq_{\Theta_n},\]
where each $\clq_{\Theta_i} = H^2(\mathbb{D})/\Theta_i
H^2(\mathbb{D})$, a Jordan block of $H^2(\mathbb{D})$ for some inner
function $\Theta_i \in H^\infty(\mathbb{D})$, or $\clq_{\Theta_i} =
H^2(\mathbb{D})$ for all $i = 1,\ldots, n$. \NI Moreover, there
exists an integer $m \in \{1, \ldots, n\}$ and inner functions
$\Theta_{i_j} \in H^{\infty}(\mathbb{D})$ such that
\[\clq^{\perp} = \mathop{\sum}_{1 \leq i_1 < \ldots < i_m \leq n}
\tilde{\Theta}_{i_j} H^2(\mathbb{D}^n),\]where
$\tilde{\Theta}_{i_j}(\bm{z})= \Theta_{i_j}(z_{i_j})$ for all
$\bm{z} \in \mathbb{D}^n$. Finally, \[P_{\clq} =
I_{H^2(\mathbb{D}^n)} - \mathop{\Pi}_{j=1}^m (I_{H^2(\mathbb{D}^n)}
- M_{\tilde{\Theta}_{i_j}} M_{\tilde{\Theta}_{i_j}}^*), \quad
\mbox{and} \quad P_{\clq^\perp} = \mathop{\Pi}_{j=1}^m
(I_{H^2(\mathbb{D}^n)} - M_{\tilde{\Theta}_{i_j}}
M_{\tilde{\Theta}_{i_j}}^*).\]
\end{Theorem}

In what follows, we realize a doubly commuting quotient module
$\clq$ of $H^2(\mathbb{D}^n)$ as $\clq_{\Theta_1} \otimes \cdots
\otimes \clq_{\Theta_n}$ where each $\clq_{\Theta_i}$ ($1 \leq i
\leq n$) is either a Jordan block of $H^2(\mathbb{D})$ (see
\cite{Ber}, \cite{B}, \cite{NF}) or the Hardy module
$H^2(\mathbb{D})$. Consequently, a co-doubly commuting submodule
$\cls$ of $H^2(\mathbb{D}^n)$ will be realized as \[\cls =
\mathop{\sum}_{1 \leq i \leq n} \tilde{\Theta}_{i}
H^2(\mathbb{D}^n),\]where $\tilde{\Theta}_{i}(\bm{z})=
\Theta_{i}(z_{i})$ for all $\bm{z} \in \mathbb{D}^n$ and each
$\Theta_{i} \in H^\infty(\mathbb{D})$ is either inner or the zero
function.

Note that a Jordan block $\clq_{\Theta}$ of $H^2(\mathbb{D})$ is of
finite dimension if and only if the inner function $\Theta$ is a
finite Blaschke products on the unit disk. Moreover, for any Jordan
block $\clq_{\Theta}$ of $H^2(\mathbb{D})$ we have
\[\mbox{rank} [C_z^*, C_z] \leq 1,\]where $C_z = P_{\clq_{\Theta}}
M_z|_{\clq_{\Theta}}$.

First, we record a simple observation concerning essentially normal
doubly commuting quotient modules of the Hardy module
$H^2(\mathbb{D}^n)$.

\begin{Proposition}\label{comm-n} Let $\clq = \clq_{\Theta_1} \otimes \cdots \otimes \clq_{\Theta_n}$ be a doubly commuting quotient
module of $H^2(\mathbb{D}^n)$. Then $\clq$ is essentially normal if
and only if each representing function $\Theta_i$ of $\clq$ is a
finite Blaschke product for all $i = 1, \ldots, n$, or,
equivalently, $\mbox{dim~} \clq < \infty$.
\end{Proposition}

\NI \textsf{Proof.} Suppose $\clq$ is a doubly commuting quotient
module of $H^2(\mathbb{D}^n)$, that is, \[[C_{z_i}^*, C_{z_j}] =
0,\] and \[C_{z_i} = I_{\clq_{\Theta_1}} \otimes \cdots \otimes
P_{\clq_{\Theta_i}} M_z|_{\clq_{\Theta_i}} \otimes \cdots \otimes
I_{\clq_{\Theta_n}},\]for all $1 \leq i < j \leq n$. Then we obtain
readily that
\[[C_{z_i}^*, C_{z_i}] = I_{\clq_{\Theta_1}} \otimes \cdots \otimes
\underbrace{[C_z^*, C_z]}_{i^{\bm\, th}} \otimes \cdots \otimes
I_{\clq_{\Theta_n}},\]and conclude that $[C_{z_i}^*, C_{z_i}] \in
\clk(\clq)$ for all $1 \leq i \leq n$ if and only if $\mbox{dim~}
\clq_{\Theta_i} < \infty$, or, equivalently, if and only if
$\Theta_i$ is a finite Blaschke product for all $i = 1, \ldots, n$.
Hence, $\clq$ is essentially normal if and only if $\Theta_i$ is a
finite Blaschke product for all $i = 1, \ldots, n$. This concludes
the proof. \qed

Hence it follows in particular that essential normality of
submodules of $H^2(\mathbb{D}^n)$ seems like a rather strong
property. Therefore, in the rest of this section we will focus only
on essentially doubly commuting submodules of $H^2(\mathbb{D}^n)$.

Before proceeding we need to prove the following result concerning
the rank of the multiplication operator restricted to a submodule
and projected back onto the corresponding quotient module of the
Hardy module $H^2(\mathbb{D})$.

\begin{Proposition}\label{1-commutator}
Let $\clq_{\Theta}$ be a quotient module of $H^2(\mathbb{D})$ for
some inner function $\Theta \in H^{\infty}(\mathbb{D})$. Then
$C_{\Theta} := P_{\clq_{\Theta}} M_z^*|_{\Theta H^2(\mathbb{D})} \in
\cll(\Theta H^2(\mathbb{D}), \clq_{\Theta})$ is given by \[C_\Theta
= [M_z^*, M_\Theta] M_{\Theta}^*.\]Moreover, $C_\Theta$ is a rank
one operator and
\[\|C_{\Theta}\| = (1 - |\Theta(0)|^2)^{\frac{1}{2}}.\]
\end{Proposition}

\NI\textsf{Proof.} We begin by calculating \[\begin{split}(I -
M_{\Theta} M_{\Theta}^*) M_z^* M_{\Theta} & = M_z^* M_{\Theta} -
M_{\Theta} M_{\Theta}^* M_z^* M_{\Theta}\\ &  = M_z^* M_{\Theta} -
M_{\Theta} M_z^*
M_{\Theta}M_{\Theta}^* \\ &= M_z^* M_{\Theta} - M_{\Theta} M_z^*\\
& = [( M_z^* M_{\Theta} - M_{\Theta} M_z^*) M_{\Theta}^*]
M_{\Theta}.\end{split}\]Therefore, we have
\[C_{\Theta} = (M_z^* M_{\Theta} - M_{\Theta} M_z^*) M_{\Theta}^* = [M_z^*,
M_\Theta] M_{\Theta}^*.\]Now for all $l \geq 1$ \[[M_z^*, M_\Theta]
z^l = (M_z^* M_\Theta - M_\Theta M_z^*) M_z^l 1 = 0,\]and \[[M_z^*,
M_\Theta] 1 = (M_z^* M_\Theta - M_\Theta M_z^*) 1 = M_z^*
\Theta.\]And so,
\[[M_z^*, M_\Theta] f = M_z^* M_{\Theta} f(0) = f(0) M_z^* \Theta =
\langle f, 1\rangle M_z^* \Theta = \langle \Theta f, \Theta \rangle
M_z^* \Theta,\]for all $f \in H^2(\mathbb{D})$. Hence, we infer that
\[C_{\Theta} (\Theta f) = [M_z^*, M_\Theta] f  =  \langle
\Theta f, \Theta \rangle M_z^* \Theta.\]Therefore, $C_\Theta$ is a
rank one operator and \[C_{\Theta} f = \langle f, \Theta \rangle
M_z^* \Theta,\]for all $f \in \Theta H^2(\mathbb{D})$. Finally,
\[\begin{split}\|C_\Theta\|^2 & = \|\Theta\|^2 \|M_z^* \Theta\|^2 = \|M_z^*
\Theta\|^2 \\& = \langle M_z M_z^* \Theta, \Theta\rangle =
\langle(I_{H^2(\mathbb{D})} - P_{\mathbb{C}}) \Theta, \Theta \rangle
\\ & = \|\Theta\|^2 - |\Theta(0)|^2
\\& = 1 - |\Theta(0)|^2.\end{split}\]This completes the proof. \qed

In the sequel we will need the following well known fact (cf. Lemma
2.5 in \cite{JS}).

\begin{Lemma}\label{P-F} Let $\{P_i\}_{i=1}^n$ be a collection of commuting orthogonal
projections on a Hilbert space $\clh$. Then \[\cll :=
\mathop{\sum}_{i=1}^n \mbox{ran} P_i,\] is closed and the orthogonal
projection of $\clh$ onto $\cll$ is given by
\[\begin{split}P_{\cll} & = P_1 (I - P_2) \cdots (I - P_n) +
P_2 (I - P_3) \cdots (I - P_n) + \cdots + P_{n-1} (I - P_n) + P_n\\
& = P_n (I - P_{n-1}) \cdots (I - P_1) + P_{n-1} (I - P_{n-2})
\cdots (I - P_1) + \cdots + P_2 (I - P_1) +
P_1.\end{split}\]Moreover,
\[P_{\cll}
= I - \mathop{\prod}_{i=1}^n (I - P_i).\]
\end{Lemma}

We now are ready to compute the cross commutators of a co-doubly
commuting submodule of $H^2(\mathbb{D}^n)$.

\begin{Theorem}\label{2-commutator}
Let $\cls = \sum_{i=1}^n \tilde{\Theta}_i H^2(\mathbb{D}^n)$ be a
co-doubly commuting submodule of $H^2(\mathbb{D}^n)$, where
$\tilde{\Theta}_i(\bm{z}) = \Theta_i(z_i)$ for all $\bm{z} \in
\mathbb{D}^n$ and each $\Theta_i \in H^{\infty}(\mathbb{D})$ is
either an inner function or the zero function and $1 \leq i \leq n$.
Then for all $1 \leq i < j \leq n$,
\[[R_{z_i}^* , R_{z_j}] = I_{\clq_{{\Theta}_1}} \otimes \cdots
\otimes \underbrace{P_{\clq_{{\Theta}_i}} M_{z}^*|_{{\Theta}_i
H^2(\mathbb{D})}}_{i^{\rm th}} \otimes \cdots \otimes
\underbrace{P_{{\Theta}_j H^2(\mathbb{D})}
M_{z}|_{\clq_{{\Theta}_j}}}_{j^{\rm th}} \otimes \cdots \otimes
I_{\clq_{{\Theta}_n}},\] and \[\|[R_{z_i}^* , R_{z_j}]\| = (1 -
|{\Theta}_i(0)|^2)^{\frac{1}{2}} (1 -
|{\Theta}_j(0)|^2)^{\frac{1}{2}}.\]
\end{Theorem}

\NI \textsf{Proof.} Let $\cls = \sum_{i=1}^n \tilde{\Theta}_i
H^2(\mathbb{D}^n)$, for some one variable inner functions $\Theta_i
\in H^{\infty}(\mathbb{D})$. Let $\tilde{P}_i$ be the orthogonal
projection in $\cll(\cls)$ defined by
\[\tilde{P}_i = M_{\tilde{\Theta}_i} M_{\tilde{\Theta}_i}^*,\] for
all $i = 1, \ldots, n$. Then it follows that
$\{\tilde{P}_i\}_{i=1}^n$ is a collection of commuting orthogonal
projections. By virtue of Theorem \ref{S-proj} and Lemma \ref{P-F},
\[\begin{split} {P}_{\cls} & = I_{H^2(\mathbb{D}^n)} - \mathop{\Pi}_{i=1}^n
(I_{H^2(\mathbb{D}^n)} - \tilde{P}_i) \\ & = \tilde{P}_1 (I -
\tilde{P}_2) \cdots (I - \tilde{P}_n) + \tilde{P}_2 (I -
\tilde{P}_3) \cdots (I - \tilde{P}_n) + \cdots + \tilde{P}_{n-1} (I
- \tilde{P}_n) + \tilde{P}_n\\ & = \tilde{P}_n (I - \tilde{P}_{n-1})
\cdots (I - \tilde{P}_1) + \tilde{P}_{n-1} (I - \tilde{P}_{n-2})
\cdots (I - \tilde{P}_1) + \cdots + \tilde{P}_2 (I - \tilde{P}_1) +
\tilde{P}_1,
\end{split}\]and
\[{P}_{\clq} = \mathop{\Pi}_{i=1}^n (I_{H^2(\mathbb{D}^n)} - \tilde{P}_i).\] On the other hand, for all $1 \leq i < j \leq n$, we
obtain
\[\begin{split}[R_{z_i}^*, R_{z_j}] & = {P}_{\cls} M_{z_i}^* M_{z_j}|_{\cls} -
{P}_{\cls} M_{z_j} {P}_{\cls} M_{z_i}^*|_{\cls},\end{split}
\]and that
\[\begin{split}{P}_{\cls} M_{z_i}^*
M_{z_j}{P}_{\cls} - {P}_{\cls} M_{z_j} {P}_{\cls}
M_{z_i}^*{P}_{\cls} & = {P}_{\cls} M_{z_i}^* M_{z_j}{P}_{\cls} -
{P}_{\cls} M_{z_j} (I - {P}_{\clq}) M_{z_i}^*{P}_{\cls}\\ & =
{P}_{\cls} M_{z_j} {P}_{\clq} M_{z_i}^* {P}_{\cls}.
\end{split} \]Furthermore we have for all $1 \leq i < j \leq n$,  \[\begin{split} {P}_{\cls} & M_{z_j} {P}_{\clq} M_{z_i}^*
{P}_{\cls} \\ & = [\tilde{P}_n (I - \tilde{P}_{n-1}) \cdots (I -
\tilde{P}_1) + \tilde{P}_{n-1} (I - \tilde{P}_{n-2}) \cdots (I -
\tilde{P}_1) + \cdots + \tilde{P}_2 (I - \tilde{P}_1) +
\tilde{P}_1]\\& \;\; \;\; M_{z_j} [\mathop{\Pi}_{l=1}^n
(I_{H^2(\mathbb{D}^n)} - \tilde{P}_l)] M_{z_i}^*
\\& \;\;\;\; [\tilde{P}_1 (I - \tilde{P}_2) \cdots (I - \tilde{P}_n) + \tilde{P}_2 (I - \tilde{P}_3) \cdots
(I - \tilde{P}_n) + \cdots + \tilde{P}_{n-1} (I - \tilde{P}_n) +
\tilde{P}_n]\\ & = [\tilde{P}_n (I - \tilde{P}_{n-1}) \cdots (I -
\tilde{P}_1) + \tilde{P}_{n-1} (I - \tilde{P}_{n-2}) \cdots (I -
\tilde{P}_1) + \cdots + \tilde{P}_2 (I - \tilde{P}_1) +
\tilde{P}_1]\\& \;\; \;\; [\mathop{\Pi}_{l \neq j}
(I_{H^2(\mathbb{D}^n)} - \tilde{P}_l)] M_{z_j} M_{z_i}^*
[\mathop{\Pi}_{l\neq i} (I_{H^2(\mathbb{D}^n)} - \tilde{P}_l)] \\&
\;\;\;\; [\tilde{P}_1 (I - \tilde{P}_2) \cdots (I - \tilde{P}_n) +
\tilde{P}_2 (I - \tilde{P}_3) \cdots (I - \tilde{P}_n) + \cdots +
\tilde{P}_{n-1} (I - \tilde{P}_n) + \tilde{P}_n]
\\ & = [\tilde{P}_j (I - \tilde{P}_{j-1}) \cdots
(I - \tilde{P}_1)] M_{z_i}^* M_{z_j} [\tilde{P}_i (I - \tilde{P}_{i+1}) \cdots (I - \tilde{P}_n)]\\
& = [(I - \tilde{P}_1) \cdots (I - \tilde{P}_{j-1}) \tilde{P}_j]
M_{z_i}^* M_{z_j} [\tilde{P}_i (I - \tilde{P}_{i+1}) \cdots (I -
\tilde{P}_n)].
\end{split}\]
These equalities shows that \[\begin{split}[R_{z_i}^*, R_{z_j}] & =
[(I - \tilde{P}_1) \cdots (I - \tilde{P}_i) \cdots (I -
\tilde{P}_{j-1}) \tilde{P}_j] M_{z_i}^* M_{z_j} [\tilde{P}_i (I -
\tilde{P}_{i+1}) \cdots (I - \tilde{P}_j) \cdots (I -
\tilde{P}_n)]\\ & = (I - \tilde{P}_1) (I - \tilde{P}_2) \cdots
(I - \tilde{P}_{i-1}) \;((I - \tilde{P}_i) M_{z_i}^* \tilde{P}_i)\; (I - \tilde{P}_{i+1}) \cdots \\
& \;\;\;\; \cdots(I - \tilde{P}_{j-1}) \;(\tilde{P}_j M_{z_j} (I -
\tilde{P}_j))\; (I - \tilde{P}_{j+1}) \cdots
(I-\tilde{P}_n).\end{split}\] Moreover,
\[[R_{z_i}^*, R_{z_j}] = [(I - \tilde{P}_1) \cdots (I - \tilde{P}_{j-1}) \tilde{P}_j]
M_{z_i}^* M_{z_j} [(I - \tilde{P}_1) \cdots (I - \tilde{P}_{i-1})
\tilde{P}_i (I - \tilde{P}_{i+1}) \cdots (I - \tilde{P}_n)],\] and
\[[R_{z_i}^*, R_{z_j}] = [(I - \tilde{P}_1) \cdots (I - \tilde{P}_{j-1}) \tilde{P}_j (I -
\tilde{P}_{j+1}) \cdots (I - \tilde{P}_n)] M_{z_i}^* M_{z_j}
[\tilde{P}_i (I - \tilde{P}_{i+1}) \cdots (I - \tilde{P}_n)].\] We
conclude that the cross-commutator $[R_{z_i}^*, R_{z_j}]$ is a
bounded linear operator from
\[\clq_{\Theta_1} \otimes \cdots \otimes \clq_{\Theta_{i-1}} \otimes
\Theta_i H^2(\mathbb{D}) \otimes \clq_{\Theta_{i+1}} \otimes \cdots
\otimes \clq_{\Theta_{j}} \otimes \cdots \otimes \clq_{\Theta_n}
\subseteq \cls\]to
\[\clq_{\Theta_1} \otimes \cdots \otimes \clq_{\Theta_i} \otimes
\cdots \otimes \clq_{\Theta_{j-1}} \otimes \Theta_j H^2(\mathbb{D})
\otimes \clq_{\Theta_{j+1}} \otimes \cdots \otimes \clq_{\Theta_n}
\subseteq \cls,\]and
\[\begin{split}[R_{z_i}^*& , R_{z_j}] = I_{\clq_{\Theta_1}} \otimes \cdots
\otimes \underbrace{{P}_{\clq_{\Theta_i}} M_{z}^*|_{\Theta_i
H^2(\mathbb{D})}}_{i^{\rm th}} \otimes \cdots \otimes
\underbrace{{P}_{\Theta_j H^2(\mathbb{D})}
M_{z}|_{\clq_{\Theta_j}}}_{j^{\rm th}} \otimes \cdots \otimes
I_{\clq_{\Theta_n}}. \end{split}\] Further, we note that  \[
\begin{split} \|[R_{z_i}^* , R_{z_j}]\| & = \|I_{\clq_{\Theta_1}}
\otimes \cdots \otimes {P}_{\clq_{\Theta_i}} M_{z}^*|_{\Theta_i
H^2(\mathbb{D})} \otimes \cdots \otimes {P}_{\Theta_j
H^2(\mathbb{D})} M_{z}|_{\clq_{\Theta_j}} \otimes \cdots \otimes
I_{\clq_{\Theta_n}}\| \\& = \|{P}_{\clq_{\Theta_i}}
M_{z}^*|_{\Theta_i H^2(\mathbb{D})}\| \|{P}_{\Theta_j
H^2(\mathbb{D})} M_{z}|_{\clq_{\Theta_j}}\|, \end{split}\] and
consequently by Proposition \ref{1-commutator} we have \[
\|[R_{z_i}^* , R_{z_j}]\| = (1 - |\Theta_i(0)|^2)^{\frac{1}{2}} (1 -
|\Theta_j(0)|^2)^{\frac{1}{2}}.\]This completes the proof. \qed

In the following corollary we reveal the significance of the
identity operators in the cross commutators of the co-doubly
commuting submodules of $H^2(\mathbb{D}^n)$ for $n \geq 2$.

\begin{Corollary}\label{EN-2}
Let $\cls = \sum_{i=1}^n \tilde{\Theta}_i H^2(\mathbb{D}^n)$ be a
submodule of $H^2(\mathbb{D}^n)$ for some one variable inner
functions $\{\tilde{\Theta}_i\}_{i=1}^n \subseteq
H^{\infty}(\mathbb{D}^n)$. Then

\NI (1) for $n = 2$: the rank of the cross commutator of $\cls$ is
at most one and the Hilbert-Schmidt norm of the cross commutator is
given by
\[\|[R_{z_1}^* , R_{z_2}]\|_{\bm\, HS} = (1 -
|{\Theta}_1(0)|^2)^{\frac{1}{2}} (1 -
|{\Theta}_2(0)|^2)^{\frac{1}{2}}.\]In particular, $\cls$ is
essentially doubly commuting.

\NI (2) for $n >2$: $\cls$ is essentially doubly commuting (or of
Hilbert-Schmidt cross-commutators) if and only if $\tilde{\Theta}_i$
is a one variable finite Blaschke product for all $1 \leq i \leq n$,
if and only if that $\cls$ is of finite co-dimension, that is,
\[\mbox{dim~} [H^2(\mathbb{D}^n)/\cls] < \infty.\]Moreover, in this case, for all $1 \leq i <
j \leq n$ \[ \|[R_{z_i}^* , R_{z_j}]\|_{\bm\, HS} = (1 -
|\Theta_i(0)|^2)^{\frac{1}{2}} (1 -
|\Theta_j(0)|^2)^{\frac{1}{2}}.\]
\end{Corollary}

Part (1) of the above corollary was obtained by R. Yang (Corollary
1.1, \cite{Y5}). We refer the reader to \cite{AC} for more details
on finite co-dimensional submodules of the Hardy modules over
$\mathbb{D}^n$.

As another consequence of the above theorem, we have the following.

\begin{Corollary}\label{<3}
Let $n > 2$ and $\cls = \sum_{i=1}^k \tilde{\Theta}_i
H^2(\mathbb{D}^n)$ be a co-doubly commuting proper submodule of
$H^2(\mathbb{D}^n)$ for some inner functions $\{\Theta_i\}_{i=1}^k$
and $k < n$. Then $\cls$ is not essentially doubly commuting.
\end{Corollary}

Combining Corollary \ref{EN-2} and Proposition \ref{comm-n} we
obtain:

\begin{Corollary}\label{pre-regid} Let
$\cls$ be a co-doubly commuting submodule of $H^2(\mathbb{D}^n)$ and
$\clq: = H^2(\mathbb{D}^n)/\cls$ and $n > 2$. Then the following are
equivalent:

(i)  $\cls$ is essentially doubly commuting.

(ii) $\cls$ is of finite co-dimension.

(iii) $\clq$ is essentially normal.
\end{Corollary}

We conclude this section with a "rigidity" result.

\begin{Corollary}\label{n=2}
Let $n  \geq 2$ and $\cls = \sum_{i=1}^n \tilde{\Theta}_i
H^2(\mathbb{D}^n)$ be an essentially normal co-doubly commuting
submodule of $H^2(\mathbb{D}^n)$ for some one variable inner
functions $\{\Theta_i\}_{i=1}^n$. If $\cls$ is of infinite
co-dimensional, then $n = 2$.
\end{Corollary}
\NI\textsf{Proof.}  The result follows from the implication (i)
$\implies$ (ii) of Corollary \ref{pre-regid}. \qed

\newsection{Representing Inner functions of Submodules}

In this section, we will obtain the explicit representations of the
Beurling-Lax-Halmos inner functions of a class of submodules of
$H^2_{\cle}(\mathbb{D})$.

Recall that a non-trivial closed subspace $\cls$ of
$H^2_{\cle}(\mathbb{D})$ is a submodule of $H^2_{\cle}(\mathbb{D})$
if and only if \[\cls = \Theta H^2_{\cle_*}(\mathbb{D}),\]for some
closed subspace $\cle_*$ of $\cle$ and inner function $\Theta \in
H^{\infty}_{\cll(\cle_*, \cle)}(\mathbb{D})$ (unique up to unitary
equivalence). This fact is known as the Beurling-Lax-Halmos theorem
and that $\Theta$ as the representing inner function of the
submodule $\cls$. Given a submodule $\cls$ of
$H^2_{\cle}(\mathbb{D})$, it is a question of interest to determine
the inner function $\Theta$ associated with $\cls$.

Now let $\cls$ be a co-doubly commuting submodule of
$H^2(\mathbb{D}^n)$. Then by Theorem \ref{S-proj} we have
\[\cls = \sum_{i=1}^n \tilde{\Theta}_i H^2(\mathbb{D}^n),\]where
$\tilde{\Theta}_i \in H^{\infty}(\mathbb{D}^n)$ is either the zero
function or one variable inner function and $i = 1, \ldots, n$. We
realize $\cls$ as a submodule of $H^2_{\cle}(\mathbb{D})$ where
$\cle = H^2(\mathbb{D}^{n-1})$. Then by the Beurling-Lax-Halmos
theorem, there exists an inner function $\Theta \in
H^{\infty}_{\cll(\cle_*, H^2(\mathbb{D}^{n-1}))}(\mathbb{D})$, for
some closed subspace $\cle_*$ of $H^2(\mathbb{D}^{n-1})$ such that
\[\cls = \sum_{i=1}^n \tilde{\Theta}_i H^2(\mathbb{D}^n) = \Theta
H^2_{\cle_*}(\mathbb{D}).\] Since
\[R_z R_z^* = M_z P_{\cls} M_z^*|_{\cls} = M_z P_{\cls} M_z^*,\]that $R_z
R_z^*$ is an orthogonal projection onto $z \cls$ and hence we have
the orthogonal projection \[P_{\cls} - R_z R_z^* = P_{\cls \ominus z
\cls}.\]On the other hand \[\begin{split}P_{\cls} - R_z R_z^* & =
M_{\Theta} M_{\Theta}^* - M_z M_{\Theta} M_{\Theta}^* M_z^* =
M_{\Theta} (I_{H^2_{\cle_*}(\mathbb{D})} - M_z M_z^*) M_{\Theta}^* =
M_{\Theta} P_{\cle_*} M_{\Theta}^*\\ & = (M_{\Theta}
P_{\cle_*})(M_{\Theta} P_{\cle_*})^*,\end{split}\]and hence
\[\begin{split}\cls \ominus z \cls & = \mbox{ran} (P_{\cls} - R_z
R_z^*) = \mbox{ran} (M_{\Theta} P_{\cle_*} M_{\Theta}^*) =
\mbox{ran} (M_{\Theta} P_{\cle_*})\\ & = \{\Theta \eta : \eta \in
\cle_*\}.\end{split}\]Note also that $\cls \ominus z \cls$ is the
wandering subspace of $\cls$, that is,
\begin{equation}\label{wandering}\cls = \overline{\mbox{span}} \{z^l (\cls
\ominus z \cls) : l \geq 0\} = \Theta
H^2_{\cle_*}(\mathbb{D}).\end{equation}

After these preliminaries we can turn to the proof of the main
result of this section.

\begin{Theorem}\label{BLH}
Let $\cls$ be a submodule of
$H^2_{H^2(\mathbb{D}^{n-1})}(\mathbb{D})$ with the
Beurling-Lax-Halmos representation $\cls = \Theta
H^2_{\cle}(\mathbb{D})$ for some closed subspace $\cle$ of
$H^2(\mathbb{D}^{n-1})$ and inner function $\Theta \in
H^{\infty}_{\clb(\cle, H^2(\mathbb{D}^{n-1}))}(\mathbb{D})$. Then
$\cls \subseteq H^2(\mathbb{D}^n)$ is a co-doubly commuting
submodule of $H^2(\mathbb{D}^n)$ if and only if there exits an
integer $m \leq n$ and orthogonal projections $\{P_2, \ldots, P_m\}$
in $\cll(H^2(\mathbb{D}))$ and an inner function $\Theta_1 \in
H^{\infty}(\mathbb{D})$ such that $\cle = H^2(\mathbb{D}^{n-1})$ and
\[\Theta(z) = {\Theta}_1(z) (I - \tilde{P}_2) \cdots (I - \tilde{P}_m) + \tilde{P}_2
(I - \tilde{P}_3) \cdots (I - \tilde{P}_m) + \cdots +
\tilde{P}_{m-1} (I - \tilde{P}_m) + \tilde{P}_m,\]for all $z \in
\mathbb{D}$, where
\[\tilde{P}_i = I_{H^2(\mathbb{D})} \otimes \cdots \otimes
\mathop{P_i}_{(i - 1)^{\bm \,th}} \otimes \cdots \otimes
I_{H^2(\mathbb{D})} \in \cll(H^2(\mathbb{D}^{n-1}).\]
\end{Theorem}

\NI\textsf{Proof.} Let $\cls$ be a co-doubly commuting submodule of
$H^2(\mathbb{D}^n)$ so that \[\cls = \mathop{\sum}_{1 \leq i_1 <
\cdots < i_m \leq n} \tilde{\Theta}_{i_j} H^2(\mathbb{D}^n),\]for
some one variable inner function $\tilde{\Theta}_{i_j} \in
H^{\infty}(\mathbb{D}^n)$ and $1 \leq i_1 < \cdots < i_m \leq n$.
Without loss of generality, we assume that $i_j = j$ for all $j = 1,
\ldots, m$, that is, \[\cls = \mathop{\sum}_{j=1}^m
\tilde{\Theta}_{j} H^2(\mathbb{D}^n).\]Then Theorem \ref{S-proj}
implies that
\[\begin{split} P_{\cls} & = I_{H^2(\mathbb{D}^n)} -
\mathop{\Pi}_{j = 1}^m (I_{H^2(\mathbb{D}^n)} - M_{\tilde{\Theta}_{j}} M_{\tilde{\Theta}_{j}}^*) \\
& = I_{H^2(\mathbb{D}^n)} - (I_{H^2(\mathbb{D}^n)} -
M_{\tilde{\Theta}_{1}} M_{\tilde{\Theta}_{1}}^*)
\mathop{\Pi}_{j=2}^m (I_{H^2(\mathbb{D}^n)} - I_{H^2(\mathbb{D})}
\otimes \tilde{P}_{j}),\end{split}\]where \[\tilde{P}_j =
\underbrace{I_{H^2(\mathbb{D})} \otimes \cdots \otimes
\mathop{M_{\Theta_j} M_{{\Theta}_j}^*}_{{(j-1)}^{\bm\,th}} \otimes
\cdots I_{H^2(\mathbb{D})}}_{(n-1) \bm \;times} \in
\cll(H^2(\mathbb{D}^{n-1})),\] for all $j = 2, \ldots, m$. Define
$\Theta \in H^{\infty}_{\cll(H^2(\mathbb{D}^{n-1}))}(\mathbb{D})$ by
\[\Theta(z) = \Theta_1(z) (I - \tilde{P}_2) \cdots (I - \tilde{P}_m) + \tilde{P}_2
(I - \tilde{P}_3) \cdots (I - \tilde{P}_m) + \cdots +
\tilde{P}_{m-1} (I - \tilde{P}_m) + \tilde{P}_m,\]for all $z \in
\mathbb{D}$. First, note that
\[M_\Theta = M_{{\Theta}_1} (I - \tilde{P}_2) \cdots (I - \tilde{P}_m) + \tilde{P}_2
(I - \tilde{P}_3) \cdots (I - \tilde{P}_m) + \cdots +
\tilde{P}_{m-1} (I - \tilde{P}_m) + \tilde{P}_m.\]Since the terms in
the sum are orthogonal projection with orthogonal ranges, we compute
\[\begin{split} M_{\Theta}^* M_{\Theta} & = M_{{\Theta}_1}^* M_{{\Theta}_1}
(I - \tilde{P}_2) \cdots (I - \tilde{P}_m) + \tilde{P}_2 (I - \tilde{P}_3) \cdots (I - \tilde{P}_m) + \cdots +
\tilde{P}_{m-1} (I - \tilde{P}_m) + \tilde{P}_m\\ & = (I -
\tilde{P}_2) \cdots (I - \tilde{P}_m) + \tilde{P}_2 (I -
\tilde{P}_3) \cdots (I - \tilde{P}_m) + \cdots + \tilde{P}_{m-1} (I
- \tilde{P}_m) + \tilde{P}_m \\& = \mathop{\Pi}_{j = 2}^m
(I_{H^2(\mathbb{D}^n)} - \tilde{P}_{j}) + (I_{H^2(\mathbb{D}^n)} -
\mathop{\Pi}_{j = 2}^m (I_{H^2(\mathbb{D}^n)} - \tilde{P}_{j})\\ & =
I_{H^2(\mathbb{D}^{n-1})},\end{split}\]and hence that $\Theta$ is an
inner function. To prove that $\Theta$ is the Beurling-Lax-Halmos
representing inner function of $\cls$, by virtue of
(\ref{wandering}), it is enough to show that
\[\overline{\mbox{span}} \{z^l \Theta H^2(\mathbb{D}^{n-1}) : l \geq
0\} = \sum_{j=1}^m \tilde{\Theta}_j H^2(\mathbb{D}^n).\]Observe that
\[\begin{split} \Theta H^2(\mathbb{D}^{n-1}) = & {\Theta}_1 (\clq_{\Theta_2}
\otimes \cdots \otimes \clq_{\Theta_m} \otimes
\underbrace{H^2(\mathbb{D}) \otimes \cdots \otimes H^2(\mathbb{D})}_{(n-m)\rm \,times}) \\
& \oplus (\Theta_2 H^2(\mathbb{D}) \otimes \clq_{\Theta_3} \otimes
\cdots \otimes
\clq_{\Theta_m} \otimes \underbrace{H^2(\mathbb{D}) \otimes  \cdots \otimes H^2(\mathbb{D})}_{(n-m)\rm\, times})\\
& \oplus \cdots \oplus (H^2(\mathbb{D}) \otimes \cdots \otimes
H^2(\mathbb{D}) \otimes \Theta_m H^2(\mathbb{D}) \otimes
\underbrace{H^2(\mathbb{D}) \otimes \cdots\otimes
H^2(\mathbb{D})}_{(n-m)\rm\, times}),\end{split}\]and hence

\[\begin{split} \overline{\mbox{span}} \{z^l & \Theta H^2(\mathbb{D}^{n-1}) : l \geq
0\} = (\Theta_1 H^2(\mathbb{D}) \otimes \clq_{\Theta_2} \otimes
\cdots \otimes \clq_{\Theta_m} \otimes \underbrace{H^2(\mathbb{D}) \otimes  \cdots \otimes H^2(\mathbb{D})}_{(n-m)\rm \,times})\\
& \;\;\; \oplus (H^2(\mathbb{D}) \otimes \Theta_2 H^2(\mathbb{D})
\otimes
\clq_{\Theta_3} \otimes \cdots \otimes \clq_{\Theta_m} \otimes \underbrace{H^2(\mathbb{D}) \otimes \cdots \otimes H^2(\mathbb{D})}_{(n-m)\rm \,times}) \oplus \\
&\;\;\;\cdots \oplus ( H^2(\mathbb{D}) \otimes \cdots \otimes
H^2(\mathbb{D}) \otimes \Theta_m H^2(\mathbb{D}) \otimes \underbrace{H^2(\mathbb{D}) \otimes \cdots \otimes H^2(\mathbb{D})}_{(n-m)\rm \,times})\\
&= \mbox{ran~} [I_{H^2(\mathbb{D}^n)} - \mathop{\Pi}_{j = 1}^m
(I_{H^2(\mathbb{D}^n)} - M_{\tilde{\Theta}_{j}}
M_{\tilde{\Theta}_{j}}^*)] \\& = \sum_{i=1}^m \tilde{\Theta}_i
H^2(\mathbb{D}^n).\end{split}\]

\NI Conversely, if $\Theta$ is given as above, then we realize
$\Theta \in H^{\infty}_{\cll(H^2(\mathbb{D}^{n-1}))}(\mathbb{D})$ by
$\tilde{\Theta} \in H^{\infty}(\mathbb{D}^n)$ where

\[\tilde{\Theta} (\bm{z}) = {\Theta_1}(z_1) (I - \tilde{P}_2) \cdots (I - \tilde{P}_m) + \tilde{P}_2
(I - \tilde{P}_3) \cdots (I - \tilde{P}_m) + \cdots +
\tilde{P}_{m-1} (I - \tilde{P}_m) + \tilde{P}_m,\]for all $\bm{z}
\in \mathbb{D}^n$. Thus
\[\tilde{\Theta} (\bm{z}) = \tilde{\Theta}_1(\bm{z}) (I -
\tilde{P}_2) \cdots (I - \tilde{P}_m) + \tilde{P}_2 (I -
\tilde{P}_3) \cdots (I - \tilde{P}_m) + \cdots + \tilde{P}_{m-1} (I
- \tilde{P}_m) + \tilde{P}_m,\]where $\tilde{\Theta}_1(\bm{z}) =
\Theta_1(z_1)$ for all $\bm{z} \in \mathbb{D}^n$. We therefore have
\[M_{\tilde{\Theta}} M_{\tilde{\Theta}}^* = \tilde{P}_1 (I - \tilde{P}_2) \cdots (I - \tilde{P}_m) + \tilde{P}_2
(I - \tilde{P}_3) \cdots (I - \tilde{P}_m) + \cdots +
\tilde{P}_{m-1} (I - \tilde{P}_m) + \tilde{P}_m,\]where $\tilde{P}_1
= M_{\tilde{\Theta}_1} M_{\tilde{\Theta}_1}^*$. Consequently,
\[M_{\tilde{\Theta}} M_{\tilde{\Theta}}^* = I_{H^2(\mathbb{D}^n)} -
\mathop{\Pi}_{i=1}^m (I_{H^2(\mathbb{D}^n)} - \tilde{P}_i),\]and
hence
\[I_{H^2(\mathbb{D}^n)} - M_{\tilde{\Theta}} M_{\tilde{\Theta}}^* =
\mathop{\Pi}_{i=1}^m (I_{H^2(\mathbb{D}^n)} -
\tilde{P}_i).\]Therefore, we conclude that\[(\mbox{ran}
M_\Theta)^{\perp} = \underbrace{(\Theta_1 H^2(\mathbb{D}))^{\perp}
\otimes (P_2 H^2(\mathbb{D}))^{\perp} \otimes \cdots \otimes (P_m
H^2(\mathbb{D}))^{\perp}}_{m \bm \; times} \otimes
\underbrace{H^2(\mathbb{D}) \otimes \cdots \otimes
H^2(\mathbb{D})}_{(n-m) \bm\; times}.\]Combine this with the
assumption that $(\mbox{ran} M_\Theta)^{\perp}$ is a quotient module
of $H^2(\mathbb{D}^n)$ to conclude that $(\mbox{ran}
M_\Theta)^{\perp}$ is a doubly commuting quotient module. This
completes the proof. \qed

The above result is a several variables generalization ($n \geq 2$)
of Theorem 3.1 in \cite{QY} by Qin and Yang.

\newsection{Rigidity of Submodules}

Let $\clm_i \subseteq H^2(\mathbb{D}^n)$, $i = 1, 2$,  be two
submodules of $H^2(\mathbb{D}^n)$. We say that $\cls_1$ and $\cls_2$
are unitarily equivalent if there exists a unitary map $U : \cls_1
\raro \cls_2$ such that \[U (M_{z_i}|_{\cls_1}) =
(M_{z_i}|_{\cls_2}) U,\]or equivalently, \[U M_{z_i} = M_{z_i}
U,\]for all $i = 1, \ldots, n$.

A consequence of Beurling's theorem ensures that, any pair of
non-zero submodules of $H^2(\mathbb{D})$ are unitarily equivalent.
The conclusion also follows directly from the unitary invariance
property of the index of the wandering subspaces associated with the
shift operators.

\NI This phenomenon is subtle, and in general not true for many
other Hilbert modules. For instance, a pair of submodules $\cls_1$
and $\cls_2$ of the Bergman modules $L^2_a(\mathbb{B}^n)$ are
unitarily equivalent if and only if $\cls_1 = \cls_2$ (see
\cite{Ri}, \cite{P}). We refer the reader to \cite{DY}, \cite{DPSY},
\cite{DF}, \cite{DK}, \cite{Y6} and \cite{G0} for more results on
the rigidity of submodules and quotient modules of Hilbert modules
over domains in $\mathbb{C}^n$.

The submodules corresponding to the doubly commuting quotient
modules also holds the rigidity property.  This is essentially a
particular case of a rigidity result due to Agrawal, Clark and
Douglas (Corollary 4 in \cite{ACD}. See also \cite{I}).

\begin{Theorem}\textsf{(Agrawal, Clark and Douglas)}
Let $\cls_1$ and $\cls_2$ be two submodules of $H^2(\mathbb{D}^n)$,
both of which contain functions independent of $z_i$ for $i = 1,
\ldots, n$. Then $\cls_1$ and $\cls_2$ are unitarily equivalent if
and only if they are equal.
\end{Theorem}

In particular, we obtain a generalization of the rigidity theorem
for $n = 2$ (see Corollary 2.3 in \cite{Y5}).

\begin{Corollary}
Let $\cls_{\Theta} = \sum_{i=1}^n \tilde{\Theta}_i
H^2(\mathbb{D})^n$ and $\cls_{\Phi} = \sum_{i=1}^n \tilde{\Phi}_i
H^2(\mathbb{D})^n$ be a pair of submodules of $H^2(\mathbb{D})^n$,
where $\tilde{\Theta}_i(\bm{z}) = \Theta_i(z_i)$ and
$\tilde{\Phi}_i(\bm{z}) = \Phi_i(z_i)$ for inner functions
$\Theta_i, \Phi_i \in H^{\infty}(\mathbb{D})$ and $\bm{z} \in
\mathbb{D}^n$ and $i = 1, \ldots, n$. Then $\cls_{\Theta}$ and
$\cls_{\Phi}$ are unitarily equivalent if and only if $\cls_{\Theta}
= \cls_{\Phi}$.
\end{Corollary}

\NI\textsf{Proof.} Clearly $\tilde{\Theta}_i \in \cls_{\Theta}$ and
$\tilde{\Phi}_i \in \cls_{\Phi}$ are independent of $\{z_1, \cdots,
z_{i-1}, z_{i+1}, \ldots, z_n\}$ for all $i = 1, \ldots, n$.
Therefore, the submodules $\cls_{\Theta}$ and $\cls_{\Phi}$ contains
functions independent of $z_i$ for all $i = 1, \ldots, n$.
Consequently, if $\cls_\Phi$ and $\cls_\Phi$ are unitarily
equivalent then $\cls_\Theta = \cls_\Phi.$ \qed

The following result is a generalization of Corollary 4.4 in
\cite{Y5} and is a consequence of the rigidity result.

\begin{Corollary}
Let $\cls_{\Theta} = \sum_{i=1}^n \tilde{\Theta}_i
H^2(\mathbb{D})^n$ be a submodules of $H^2(\mathbb{D})^n$, where
$\tilde{\Theta}_i(\bm{z}) = \Theta_i(z_i)$ for inner functions
$\Theta_i \in H^{\infty}(\mathbb{D})$ for all $i = 1, \ldots, n$ and
$\bm{z} \in \mathbb{D}^n$. Then $\cls_{\Theta}$ and
$H^2(\mathbb{D}^n)$ are not unitarily equivalent.
\end{Corollary}

\NI\textsf{Proof.} The result follows from the previous theorem
along with the observation that $\cls_{\Theta}^{\perp} \neq
\{0\}$.\;\; \qed

We close this section by noting that the results above are not true
if we drop the assumption that all $\Theta_i$ are inner. For
instance, if $\Theta_i = \Phi_i = 0$ for all $i \neq 1$ then
$\cls_\Theta \cong \cls_\Phi$ but in general, $\cls_\Theta \neq
\cls_\Phi$ (see \cite{M}).

\newsection{Concluding remarks}

One of the central issues in the study of Hilbert modules is the
problem of analyzing essentially normal submodules and quotient
modules of a given essentially normal Hilbert module over
$\mathbb{C}[\bm{z}]$. There is, however, a crucial difference
between the Hilbert modules of functions defined over the unit ball
and the polydisc in $\mathbb{C}^n$. For instance, a submodule $\cls$
of an essentially normal Hilbert module $\clh$ is essentially normal
if and only if the quotient module $\clh/\cls$ is so (see \cite{A},
\cite{D}), that is, the study of essentially normal submodules and
quotient modules of essentially normal Hilbert modules amounts to
the same. However, this is not the case for the study of essentially
doubly commuting Hilbert modules over $\mathbb{D}^n$. In other
words, the theory of essentially doubly commuting submodules and
quotient modules of an essentially doubly commuting Hilbert modules
are two different concepts.

One could, however, consider the co-doubly commuting submodules as a
special class of submodules of the Hardy module and the results of
this paper indicates that the general picture of essentially doubly
commuting submodules of the Hardy module will by no means be easy to
understand (cf. Corollary \ref{EN-2}). In particular, the homogenous
submodules of $H^2(\mathbb{D}^2)$ are always essentially doubly
commuting \cite{CMY}. Hence Question \ref{Q1} has an affirmative
answer for the class of homogenous submodules of
$H^2(\mathbb{D}^2)$. It is not known whether the homogeneous
submodules of $H^2(\mathbb{D}^n)$, when $n \geq 3$, are essentially
doubly commuting. Corollary \ref{EN-2} gives an indication of a
possible answer to the case $n \geq 3$. Results related to
essentially normal submodules of the Drury-Arveson module over the
unit ball of $\mathbb{C}^2$ can be found in \cite{G}.

Our result concerning the Beurling-Lax-Halmos inner function,
Theorem \ref{BLH}, is closely related to the classification theory
of multi-isometries (see \cite{BCL} and \cite{BDF}) for $n = 2$
case. We hope to discuss the general case in a future paper.

 We conclude with a result concerning
the $C_0$ class. Recall that a completely non-unitary contraction
$T$ on some Hilbert space $\clh$ is said to be in the class $C_0$ if
there is a non-zero function $\Theta \in H^\infty(\mathbb{D})$ such
that $\Theta(T) = 0$ \cite{NF}.

\begin{Proposition}
Let $\clq$ be a non-trivial doubly commuting quotient module of
$H^2(\mathbb{D}^n)$. Then $R_{z_i} \in C_0$ for some $1 \leq i \leq
n$.
\end{Proposition}

\NI\textsf{Proof.} By virtue of Theorem \ref{S-proj}, we let $\clq =
\clq_{\Theta_1} \otimes \cdots \otimes \clq_{\Theta_n}$ and
$\clq_{\Theta_i} \neq H^2(\mathbb{D})$ for some $1 \leq i \leq n$.
Consequently, $\Theta_i(R_z) = 0$ and hence
\[\tilde{\Theta}_i(R_{z_i}) = I_{\clq_{\Theta_1}} \otimes \cdots \otimes
\underbrace{\Theta_i(R_z)}_{i^{\bm\,th}} \otimes \cdots \otimes
I_{\clq_{\Theta_n}} = 0.\]This concludes the proof. \qed

The above result for the case $n = 2$ is due to Douglas and Yang
(see Proposition 4.1 in \cite{DY1}). However, our proof is more
elementary.

\end{document}